\documentclass[letterpaper, 10 pt, conference]{ieeeconf}  

\IEEEoverridecommandlockouts                              

\usepackage{amsmath,amssymb}
\usepackage[hidelinks]{hyperref}
\usepackage{xcolor}
\usepackage{graphicx}
\usepackage{booktabs}
\usepackage{dsfont}
\usepackage{algorithm}
\usepackage{algpseudocode}
\usepackage[style=ieee,backend=bibtex,url=false,doi=false,isbn=false,date=year,citestyle=numeric-comp,maxbibnames=10,minbibnames=10,maxcitenames=10,mincitenames=10]{biblatex}

\usepackage{xcolor} 
\usepackage{cuted}

\bibliography{paper}

\newtheorem{defi}{Definition}
\newtheorem{ass}{Assumption}
\newtheorem{thm}{Theorem}

\newtheorem{rem}{Remark}
\newtheorem{lem}{Lemma}

\newtheorem{ex}{Example}

\DeclareMathOperator*{\argmin}{argmin} 

\title{\LARGE \bf
	Decentralized non-convex optimization via bi-level SQP and ADMM\thanks{This document is a revised version of the corresponding conference paper which appeared in the proceedings of the 2022 IEEE 61st Conference on Decision and Control, pp. 273-278, doi: 10.1109/CDC51059.2022.9992379. 
	Compared to the original conference paper, we have introduced Assumption~\ref{ass:decoupled} to correct an error in Lemma~\ref{lem:admmConvergence} and have adjusted Theorem~\ref{thm:dsqpConv} accordingly.
	Moreover, this updated version strenghtens the assumption on the continuous differentiability of the components of NLP~\eqref{eq:sepForm}, discusses the uniqueness of the SQP steps in Theorem~\ref{thm:sqpConv}, and corrects a mistake in the proof of Lemma~\ref{lem:admmAS}. Finally, we have added Remark~\ref{rem:dualinit} and Example~\ref{ex:decoupling} and have made minor changes to notation to improve readability.}}

\author{Gösta Stomberg, Alexander Engelmann and Timm Faulwasser
	\thanks{The authors were with the Institute of Energy Systems, Energy Efficiency and Energy Economics, TU Dortmund University, Dortmund, Germany. GS and TF are now with the Institute of Control Systems, Hamburg University of Technology, Hamburg, Germany. AE is now with logarithmo GmbH \& Co. KG, Dortmund, Germany.
	E-Mail: goesta.stomberg@tu-dortmund.de, \{alexander.engelmann,timm.faulwasser\}@ieee.org.}
}

\begin{document}
	
	\maketitle
	\thispagestyle{empty}
	\pagestyle{empty}
	
	\begin{abstract} 
			Decentralized non-convex optimization is important in many problems of practical relevance.
			Existing decentralized methods, however, typically either lack convergence guarantees for general non-convex problems, or they suffer from a high subproblem complexity.
			We present a novel  bi-level SQP method, where the inner quadratic problems are solved via ADMM.
			A decentralized stopping criterion from inexact Newton methods allows the early termination of ADMM as an inner algorithm to improve computational efficiency.
			The method has local convergence guarantees for non-convex problems. Moreover, it only solves sequences of Quadratic Programs, whereas many existing algorithms solve sequences of Nonlinear Programs.
			The method shows competitive numerical performance for an optimal power flow problem.
	\end{abstract}
	
	\section{Introduction}

Decentralized optimization methods, i.e. methods which solve optimization problems purely based on neighbor-to-neighbor communication, are of interest in many applications such as Optimal Power Flow (OPF)  and distributed Model Predictive Control (MPC) \cite{Christakou2017,Stewart2011}.
In many cases, these applications require to solve optimization problems with non-convex objectives and non-convex constraints.
At large, existing algorithms can be categorized as follows: a) they are decentralized, but lack guarantees for problems with non-convex constraints; b) they are not decentralized, i.e. they require centralized communication/coordination; or c) they solve non-convex Nonlinear Programs (NLP)s at a subsystem level, which increases complexity and impedes convergence guarantees.

The Alternating Direction Method of Multipliers (ADMM) is a decentralized method which shows promising performance for a large variety of problems.
ADMM is guaranteed to converge only for special classes of non-convex problems~\cite{Wang2019}.
Decentralized schemes with global convergence guarantees based on augmented Lagrangian methods are presented in \cite{Hours2016,Sun2019a}.
Whereas \cite{Hours2016} assumes polynomial objectives and equality constraints, \cite{Sun2019a} allows for more general constraints.
A distributed method with local convergence guarantees for non-convex probelms is the Augmented Lagrangian Alternating Direction Inexact Newton (ALADIN) method~\cite{Houska2016}.
Bi-level ALADIN variants employ inner algorithms to decompose the coordination step of ALADIN  \cite{Engelmann2020c}.
With the exception of~\cite{Hours2016}, all the above methods require solving constrained non-convex NLPs to optimality, which can often not be  guaranteed by numerical solvers for an arbitrary initialization~\cite{Nocedal2006}.\footnote{Most solvers are  guaranteed to converge to a stationary point of a merit function, which is not necessarily a minimizer, cf. e.g.  \cite[Theorem 19.2]{Nocedal2006}.}

An essentially decentralized interior point method with local convergence guarantees for general non-convex problems is proposed in \cite{Engelmann2021a}. 
Therein, the solution of NLPs by subsystems is avoided. However, the step-size selection requires scalar global communication.
Tailored algorithms for distributed NMPC can be found e.g. in \cite{Necoara2009a,Stewart2011}.
These algorithms also suffer from either of the above drawbacks.

For centralized non-convex optimization, Sequential Quadratic Programming (SQP) 
has local or global convergence guarantees~\cite{Nocedal2006}.
Bi-level algorithms with SQP on the outer level and ADMM as an inner algorithm are presented in \cite{Liu2020b,Lu2015}. 
The scheme of \cite{Liu2020b} has global convergence guarantees for problems with non-convex objective and linearly coupled subproblems. 
However, it does not allow for subsystem constraints or inequality constraints.
In \cite{Lu2015}, non-convex objectives and non-convex local constraints are considered, but no convergence guarantees are established.
Both approaches assume that ADMM solves the subproblem Quadratic Programs (QP)s to optimality in each SQP step, which is unrealistic in many applications. 
A further bi-level approach with ADMM as an inner method for sequential convex programming is presented in~\cite{Le2020}.

The focus of this work is a bi-level decentralized SQP method (d-SQP) with ADMM as inner algorithm.
As the proposed scheme combines SQP with inner decentralized ADMM we obtain a fully decentralized algorithm.
This paper presents two contributions. 
First, in contrast to the bi-level SQP method in \cite{Liu2020b}, we establish convergence for nonlinear programs with non-convex objective \emph{and} non-convex constraints. 
Second, we derive a novel stopping criterion for the inner ADMM iterations that guarantees local convergence of d-SQP despite inexact QP solutions. 
This is crucial as it may not be computationally feasible to obtain an exact subproblem solution via ADMM in each SQP step.

The paper is organized as follows: Section~\ref{sec:problem} introduces the problem formulation, SQP schemes, and ADMM. Section~\ref{sec:dsqp} derives d-SQP and establishes local convergence.
Section~\ref{sec:numerics} analyzes numerical results from OPF. 
	
\textit{Notation:} Given a matrix $A$ and an integer $j$, $[A]_j$ denotes the $j$th row of $A$. For an index set $\mathcal{A}$, $[A]_\mathcal{A}$ denotes the matrix consisting of rows $ [A]_j$ for all $j \in \mathcal{A}$. Likewise, $[a]_j$ is the $j$th component of vector $a$ and $a_\mathcal{A}$ is the vector of components $[a]_j$ for all $j \in \mathcal{A}$. The concatenation of vectors $x$ and $y$ into a column vector is $(x,y)$. 
The symbol $\| \cdot \|$ denotes any vector norm on $\mathbb{R}^n$ or its induced matrix norm, respectively.
The open $\varepsilon$-neighborhood around a point $x$ is denoted by $\mathcal{B}_\varepsilon(x)$, i.e., ${\mathcal{B}_\varepsilon(x) \doteq \{y \in \mathbb{R}^{n_x} | \| y - x\| < \varepsilon\}}$. The condition number of a matrix $A \in \mathbb{R}^{n \times n}$ is defined as ${\text{cond}(A)\doteq \| A \| \| A^{-1} \|}$. 
$I$ is the identity matrix of appropriate dimension.
Given two vectors $a,b \in \mathbb{R}^n$, we denote the vector of componentwise minima by $\min(a,b) \doteq (\min([a]_1,[b]_1),\dots,\min([a]_n,[b]_n))$.

\section{Problem statement and preliminaries}\label{sec:problem}
We consider Nonlinear Programs of the form
\begin{subequations} \label{eq:sepForm}
	\begin{align} 
	\min_{x_1,\dots,x_S} \; \sum_{i \in \mathcal{S}} &f_i(x_i) \\
	\text{subject to}\hspace{0.5mm} \quad  g_i(x_i)&=0 \; | \; \nu_i & \forall i \in \mathcal{S}, \label{eq:sepProbGi} \\
	h_i(x_i) &\leq 0 \; | \; \mu_i & \forall i \in \mathcal{S},\label{eq:sepProbHi} \\
	\sum_{i \in \mathcal{S}} E_ix_i  &= c \hspace*{1.2mm} | \; \lambda,\label{eq:consConstr}
	\end{align}
\end{subequations} where $\mathcal{S} = \{1,\dots,S\}$ is a set of subsystems, each of~which is equipped with decision variables
$x_i \hspace*{-0.1cm}\in \hspace*{-0.1cm}\mathbb{R}^{n_i}$ and three times continuously differentiable functions
$f_i: \mathbb{R}^{n_i} \rightarrow \mathbb{R}$,
$g_i:\mathbb{R}^{n_i} \rightarrow \mathbb{R}^{n_{g,i}}$, and
${h_i:\mathbb{R}^{n_i} \rightarrow  \mathbb{R}^{n_{h,i}}}$. 
The matrices $E_i \in \mathbb{R}^{{n_c} \times {n_i} }$ and the vector ${c \in \mathbb{R}^{n_c}}$ couple the subsystems.
The notation in \eqref{eq:sepForm} highlights that $\nu_i \in \mathbb{R}^{n_{g,i}}$, $\mu_i \in \mathbb{R}^{n_{h,i}}$, and $\lambda \in \mathbb{R}^{n_c}$ are Lagrange multipliers associated with the respective constraints.
The centralized variables are ${x\doteq (x_1,\dots,x_S)}$, $\nu\doteq (\nu_1,\dots,\nu_S)$, and $\mu\doteq (\mu_1,\dots,\mu_S)$.
We define the Lagrangian of \eqref{eq:sepForm},
\begin{align*}
L(&x,\nu,\mu,\lambda) = \left(\sum_{i \in \mathcal{S}} L_i(x_i,\nu_i,\mu_i,\lambda) \right)- \lambda^\top c,
\end{align*} where $L_i(\cdot) = f_i(x_i) + \nu_i^\top g_i(x_i) + \mu_i^\top h_i(x_i) + \lambda^\top E_i x_i$.
Then, the Karush-Kuhn-Tucker (KKT) conditions of \eqref{eq:sepForm} read
\begin{subequations}\label{eq:KKT}
	\begin{align}
	\nabla f_i(x_i) \hspace*{-0.7mm} + \hspace*{-0.7mm}\nabla g_i(x_i) \nu_i \hspace*{-0.7mm}+\hspace*{-0.7mm} \nabla h_i(x_i) \mu_i \hspace*{-0.7mm}+\hspace*{-0.7mm} E_i^\top \lambda &= 0 \;\;  \forall i \in \mathcal{S},\\
	g_i(x_i) &= 0\;\; \forall  i \in \mathcal{S},\\
	h_i(x_i) \leq 0, \; \mu_i \geq 0, \; \mu_i^\top h_i(x_i) &= 0\;\; \forall  i \in \mathcal{S},\label{eq:NLP_KKT_ineq} \\
	\sum_{i \in \mathcal{S}} E_i x_i &= c. 
	\end{align}
\end{subequations}
Throughout the paper, we denote the primal-dual variables of NLP~\eqref{eq:sepForm} by $p\doteq(x,\nu,\mu,\lambda)$ and a KKT point by ${p^\star\doteq (x^\star,\nu^\star,\mu^\star,\lambda^\star)}$.
Moroever, let $\mathcal{A}_i$ and $\mathcal{I}_i$ denote the sets of active and inactive inequality constraints at $x_i^\star$ respectively, i.e.,
\begin{align*}
\mathcal{A}_i &\doteq \{j \in \{1,\dots,n_{h,i} \} | [h_i(x_i^\star)]_j = 0\} \quad \forall i \in \mathcal{S}\\
\mathcal{I}_i &\doteq \{j \in \{1,\dots,n_{h,i} \} | [h_i(x_i^\star)]_j < 0\} \quad \forall i \in \mathcal{S}.
\end{align*}

\subsection{Sequential quadratic programming}

SQP methods repeatedly solve quadratic approximations of~\eqref{eq:sepForm}, cf. \cite[Ch. 18]{Nocedal2006}.
A quadratic approximation of  \eqref{eq:sepForm} at a primal-dual iterate $p^k = (x^k,\nu^k,\mu^k,\lambda^k)$ is
\begin{subequations}\label{eq:qpk}
	\begin{align} 
	\min_{s_1,\dots,s_S} \; \sum_{i \in \mathcal{S}} \bigg( \frac{1}{2} s_i^\top H_{i}^k  & s_i + \nabla f_i^{k\top} s_i \bigg)\\
	\nonumber\text{subject to}\hspace{1.25cm}&\\
	g_i^k + \nabla g_i^{k\top} s_i&=0 \; |\; \nu_i \quad \hspace*{0.3mm} \forall i \in \mathcal{S},\\
	h_i^k + \nabla h_i^{k\top} s_i &\leq 0 \; |\; \mu_i \quad \forall i \in \mathcal{S},\\
	\sum_{i \in \mathcal{S}} E_i(x_i^k + s_i) &= c \hspace*{1.2mm} |\; \lambda,
	\end{align}
\end{subequations} where $H_{i}^k \doteq \nabla_{x_ix_i}^2 L_i(x_i^k,\nu_i^k,\mu_i^k,\lambda^k)$.
The notation $g_i^k$ and $\nabla g_i^k$ is shorthand for $g_i(x_i^k)$ and $\nabla g_i(x_i^k)$, respectively and the same applies to the functions $f_i$ and $h_i$.
Observe that the KKT system~\eqref{eq:KKT} can be written as $F(x,\nu,\mu,\lambda) = 0$, where
\begin{equation}\label{eq:F}
F(x,\nu,\mu,\lambda) \doteq \begin{bmatrix}
\nabla_{x_1} L_1(x_1,\nu_1,\mu_1,\lambda)\\
g_1(x_1)\\
\text{min}(-h_1(x_1),\mu_1)\vspace*{-0.2cm}\\
\vdots\vspace*{-0.2cm}\\
\nabla_{x_S} L_S(x_S,\nu_S,\mu_S,\lambda)\\
g_S(x_S)\\
\text{min}(-h_S(x_S),\mu_S)\\
\left(\sum_{i \in \mathcal{S}}E_i x_i \right)- c
\end{bmatrix}.
\end{equation}
The block rows $\text{min}(-h_i(x_i),\mu_i)$ represent \eqref{eq:NLP_KKT_ineq}.
Algorithm~\ref{alg:basicSQP} summarizes an SQP method for solving \eqref{eq:sepForm}.

\begin{algorithm}[t]
	\caption{Inequality-constrained SQP for solving~\eqref{eq:sepForm}}
	\begin{algorithmic}[1]
		\State Initialization: $k = 0$, $(x_i^0,\nu_i^0,\mu_i^0)$ for all $i \in \mathcal{S}$, $\lambda^0$, $\epsilon$ \label{stp:1}
		\While{ $ \| F^k \| \nleq \epsilon $ } 
		\State compute $\nabla f_i^k, g_i^k, \nabla g_i^k, h_i^k, \nabla h_i^k, H_{i}^k\, \, \forall i \in \mathcal{S}$ 
		\State $(s^k,\nu^{k+1},\mu^{k+1},\lambda^{k+1}) \leftarrow$ solve QP \eqref{eq:qpk} \label{stp:7}
		\State $x_i^{k+1} = x_i^k  + s_i^k$ for all $i \in \mathcal{S}$
		\State $k \leftarrow k +1$
		\EndWhile
		\State \textbf{return} $x_i^k$ for all $i \in \mathcal{S}$
	\end{algorithmic} \label{alg:basicSQP}
\end{algorithm} 

\begin{ass}[Regular KKT point]\label{ass:kkt}
	The point $p^\star$ is a KKT point of \eqref{eq:sepForm} which, for all $i \in \mathcal{S}$, satisfies
	\begin{enumerate}
		\item[i)] $h_i(x_i^\star) + \mu_i^\star \neq 0$ (strict complementarity), 
		\item[ii)] $s_i^\top \nabla_{x_ix_i}^2 L_i(x_i^\star,\nu_i^\star,\mu_i^\star,\lambda^\star) s_i > 0$ for all $s_i \neq 0$ with $\nabla g_i(x_i^\star)^\top s_i = 0$.\footnote{This is a slightly stronger assumption than the Second-Order Sufficient Condition as we exclude the conditions ${[\nabla h_i(x_i^\star)]_{\mathcal{A}_i}^\top s_i = 0}$ and $E_is_i=0$.}
	\end{enumerate}
	Furthermore, the matrix
	\begin{align*}
	\begin{bmatrix} \nabla g_1(x_1^\star)^\top &  & \\
	& \ddots & \\
	& & \nabla g_S(x_S^\star)^\top\\
	[\nabla h_1(x_1^\star)^\top]_{\mathcal{A}_1} & & \\
	& \ddots & \\
	& & [\nabla h_S(x_S^\star)^\top]_{\mathcal{A}_S}\\
	E_1 & \dots & E_S \end{bmatrix}
	\end{align*} has full row rank, i.e., it satisfies the Linear Independence Constraint Qualification (LICQ).\hfill $\square$ 
\end{ass}
\begin{defi}[Convergence rates]
	We say that the sequence $\{p^k\} \subset \mathbb{R}^{n_p}$ converges to $p^\star \in \mathbb{R}^{n_p}$
	\begin{enumerate}
		\item[i)] q-linearly, if ${\| p^{k+1} - p^\star \| \leq c \| p^k - p^\star \| } \quad \forall k \geq k_0$ for some $0<c<1$ and $k_0 \geq 0$.
		\item[ii)] q-\{superlinearly, quadratically\}, if ${p^k \rightarrow p^\star}$ and $\|p^{k+1} - p^\star \| = \{ o (\| p^k - p^\star \|), O(\|p^{k} - p^\star \|^2) \} $ for $k \rightarrow \infty$.
	\end{enumerate}
\end{defi}
\begin{thm}[{Local convergence of SQP}]\label{thm:sqpConv}
 Let Assumption~\ref{ass:kkt} hold. Then, there exists a constant $\varepsilon_1 > 0$ such that, for all $p^0 \in \mathcal{B}_{\varepsilon_1} ( p^\star)$, the sequence $\{p^k\}$ generated by Algorithm~\autoref{alg:basicSQP} converges q-quadratically to $p^\star$. \hfill $\square$
\end{thm}
\begin{proof}
	We first prove that, (a), QP~\eqref{eq:qpk} has a unique solution $p^{k,\star}$ if $p^k \approx p^\star$. Then, (b), we invoke a classic result from~\cite{Geiger2002} to obtain convergence.

	(a) Assumption~\ref{ass:kkt} yields that, if $p^k = p^\star$, then the solution $p^{k,\star}$ to QP~\eqref{eq:qpk} satisfies LICQ, strict complementarity, and the stronger Second-Order Sufficient Condition~(SOSC) Assumption~\ref{ass:kkt} ii).
	We next view QP~\eqref{eq:qpk} formed at $p^k$ as a perturbed version of QP~\eqref{eq:qpk} formed at $p^\star$ and apply the Basic Sensitivity Theorem~(BST,~\cite[Theorem 3.2.2]{Fiacco1983}).
	Thus, if $p^k \in \mathcal{B}_{\varepsilon_1}(p^\star)$ for sufficiently small $\varepsilon_1 > 0$, then any solution to QP~\eqref{eq:qpk} satisfies LICQ and strict complementarity with the same active set as NLP~\eqref{eq:sepForm}, i.e.
	\begin{align*}
	[h_i^k]_j + [\nabla h_i^{k\top} s_i^{k,\star}]_j  &= 0 \quad \forall j \in \mathcal{A}_i, \; \forall i \in \mathcal{S},\\
	[\mu_i^{k,\star}]_j & = 0 \quad \forall j \in \mathcal{I}_i, \hspace*{1.8mm} \forall i \in \mathcal{S}.
	\end{align*}
	Moreover, for all $i \in \mathcal{S}$,
	\begin{equation}\label{eq:sconvQP}
	s_i^\top H_i^k s_i > 0 \quad \forall s \neq 0 \quad \text{with} \quad \nabla g_i^{k\top} s_i = 0,
	\end{equation} which follows by adjusting the proof of the BST in~\Cite{Fiacco1983} to the stronger SOSC Assumption~\ref{ass:kkt} ii).
	Therefore, the objective of QP~\eqref{eq:qpk} is strictly convex over the constraints and hence $p^{k,\star}$ is unique.

	(b) Since we assume the components $f_i$, $g_i$, and $h_i$ of NLP~\eqref{eq:sepForm} to be three times continuously differentiable for all $i \in \mathcal{S}$, we obtain that the Hessians $\nabla^2 f_i$, $\nabla^2 [g_i]_j$, $j \in \{1,\dots,n_{g,i}\}$, and $\nabla^2 [h_i]_j$, $j \in \{1,\dots,n_{h,i}\}$ are Lipschitz continuous for $p^k \approx p^\star$ and for all $i \in \mathcal{S}$.   
	Local q-quadratic convergence for sufficiently small $\varepsilon_1 > 0$ thus follows from~\cite[Theorem 5.31]{Geiger2002}, where we have inserted the uniqueness of the QP solution derived in~(a).
\end{proof}
	
The proof of Theorem~\ref{thm:sqpConv} given in \cite{Geiger2002} shows that, inside $\mathcal{B}_{\varepsilon_1} (p^\star)$, the iterates generated by Algorithm~\ref{alg:basicSQP} correspond to the iterates generated by Newton's method applied to the nonlinear system of equations ${F(p) = 0}$.
In particular, for all $p^k \in \mathcal{B}_{\varepsilon_1}(p^\star)$, $F^k\doteq F(p^k)$ is continuously differentiable and the Jacobian $\nabla F^{k\top}\doteq \nabla F(p^k)^\top$ is regular. 
Hence the Newton iteration $p^{k+1} = p^k + d^k$ with $F^k  + \nabla F^{k\top} d^k = 0$, where
\begin{equation}\label{eq:d}
d^k \doteq ( s^k, \Delta \nu^k, \Delta \mu^k, \Delta \lambda^k),
\end{equation} ${\Delta \nu^k \doteq \nu^{k+1} - \nu^k}$, $\Delta \mu^k \doteq \mu^{k+1} - \mu^k$, and $\Delta \lambda^k \doteq \lambda^{k+1} - \lambda^k$ is well-defined inside $\mathcal{B}_{\varepsilon_1}$.

\subsection{Alternating direction method of multipliers}

We next recall how to solve \eqref{eq:qpk} via ADMM. First, we reformulate \eqref{eq:qpk} as
\begin{subequations}\label{eq:reformulatedQP}
	\begin{align}
	\min_{\substack{s_1 \in \mathbb{S}_1^k,\dots, s_S \in \mathbb{S}_s^k,\\\bar{s}_1,\dots,\bar{s}_S}} \sum_{i \in \mathcal{S}} \phi_i^k(s_i)\\
	\text{subject to}\hspace{0.5mm} \quad s_i - \bar{s}_i &=0 \; |\; \gamma_i^{\phantom{QP}} \quad \forall i \in \mathcal{S}, \label{eq:consConstr2}\\
	\sum_{i \in \mathcal{S}} E_i(x_i^k + \bar{s}_i) &= c \hspace*{1.2mm} |\; \lambda,
	\end{align}
\end{subequations} where $\bar{s}_i \in \mathbb{R}^{n_{i}}$ is an auxiliary decision variable for each subsystem $i \in \mathcal{S}$, $\phi_i^k(s_i)\doteq s_i^\top H_{i}^k s_i/2 + \nabla f_i^{k\top} s_i $, and 
\begin{align*}
\mathbb{S}_i^k \doteq \left\{ s_i \in \mathbb{R}^{n_i} \; \left| \; \begin{aligned} g_i^k + \nabla g_i^{k\top} s_i &= 0\\
h_i^k + \nabla h_i^{k\top} s_i &\leq 0 \end{aligned} \right. \right\}.
\end{align*} 
Let $\bar{s} \hspace*{-0.5mm} \doteq \hspace*{-0.5mm} (\bar{s}_1,\dots,\bar{s}_S)$, ${s \hspace*{-0.5mm}\doteq \hspace*{-0.5mm} ({s}_1,\dots,s_S)}$, and $\gamma \hspace*{-0.5mm} \doteq \hspace*{-0.5mm} (\gamma_1,\dots,\gamma_S)$.
The augmented Lagrangian of \eqref{eq:reformulatedQP} with respect to \eqref{eq:consConstr2} reads
\begin{align*}
L_\rho^k(s,\bar{s},\gamma) &= \sum_{i \in \mathcal{S}}  L_{\rho,i}^k(s_i,\bar{s}_i,\gamma_i)\\
&=\sum_{i \in \mathcal{S}} \left( \phi_i^k (s_i) + \gamma_i^\top (s_i - \bar{s}_i) + \frac{\rho}{2}\| s_i - \bar{s}_i \|^2_2 \right).
\end{align*}
Denote the feasible set of the coupling constraints by ${ \mathbb{E} = \{ \bar{s} \in \mathbb{R}^n \;|\; \sum_{i \in \mathcal{S}} E_i (x_i^k+\bar{s}_i)=c\} }$.
Algorithm~\ref{alg:ADMM} summarizes ADMM, where the dual variables $\nu_i^{l+1}$ and $\mu_i^{l+1}$ are obtained in Step~\ref{eq:ADMMStep1}.
Step~\ref{eq:ADMMStep1} can be executed by each subsystem in parallel by solving, for $i \in \mathcal{S}$, the subsystem QP 
\begin{subequations} \label{eq:localQP}
	\begin{align} 
	\min_{s_i} \; \frac{1}{2} s_i^\top (H_i^k + \rho I) s_i + &(\nabla f_i^k + \gamma_i^l - \rho \bar{s}_i^l)^{\top} s_i \\
	\nonumber \textrm{subject to}\hspace{2.1cm}&\\
	\quad  g_i^k + \nabla g_i^{k\top} s_i&=0 \; |\; \nu_i,\label{eq:locEqCon}\\
	h_i^k + \nabla h_i^{k\top} s_i &\leq 0 \; | \;\mu_i,\label{eq:locIneqCon}
	\end{align}
\end{subequations} where the iterates $\bar{s}_i^l$ and $\gamma_i^l$ are parameters.

\begin{rem}[Decentralized ADMM]
Step~\ref{eq:ADMMStep2} in Algorithm~\ref{alg:ADMM} solves a QP subject to the coupling constraints and in general requires centralized computation. 
However, if \eqref{eq:sepForm} is given in so-called consensus form, and 
if~$\gamma_i^0$ 
is chosen appropriately, then Step~\ref{eq:ADMMStep2} is equivalent to a decentralized averaging step, i.e., it only requires neighbor-to-neighbor communication and local computation~\cite[Ch. 7]{Boyd2011}. 
A common setup where Problem~\eqref{eq:sepForm} is given in consensus form occurs if the constraints~\eqref{eq:consConstr} couple original and copied variables between neighboring subsystems, cf. \cite{Stomberg2022}. Then, if the SQP initialization satisfies $\sum_{i \in \mathcal{S}} E_i x_i^0 = c$, we may replace Step~\ref{eq:ADMMStep2} by a decentralized averaging step. $\hfill \square$
\end{rem}	

\begin{algorithm}[t]
	\caption{ADMM for solving \eqref{eq:reformulatedQP}}
	\begin{algorithmic}[1]
		\State Initialization: $l = 0$, $(\bar{s}_i^0,\gamma_i^0)$ for all $i \in \mathcal{S}$ 
		\While{ not converged}
		\State $\displaystyle (s_i^{l+1},\nu_i^{l+1},\mu_i^{l+1}) \leftarrow \min_{s_i \in \mathbb{S}_i^k} L_{\rho,i}^k(s_i,\bar{s}_i^l,\gamma_i^l)\,\forall i \in \mathcal{S}$ \label{eq:ADMMStep1}  
		\State $\displaystyle \bar{s}^{l+1} = \argmin_{\bar{s} \in \mathbb{E}} \sum_{i \in \mathcal{S}} L_{\rho,i}^k(s_i^{l+1},\bar{s}_i,\gamma_i^{l})$ \label{eq:ADMMStep2}
		\State $\gamma_i^{l+1} = \gamma_i^l + \rho (s_i^{l+1}-\bar{s}_i^{l+1})$ for all $i \in \mathcal{S}$ \label{eq:ADMMStep3}
		\State $l \leftarrow l+1$
		\EndWhile
		\State \textbf{return} $(\bar{s}_i^l,\nu_i^{l},\mu_i^{l},\gamma_i^l)$ \textbf{for all} $i \in \mathcal{S}$
	\end{algorithmic} \label{alg:ADMM}
\end{algorithm}
	
\section{Decentralized sequential quadratic programming}\label{sec:dsqp}
The key idea of our  approach is to solve QP~\eqref{eq:qpk} with ADMM. 
It may not be tractable to solve \eqref{eq:qpk} to high precision in every SQP step.
Therefore we next derive convergence results of the inexact outer SQP steps and then provide technical results for the subproblems based on well-known ADMM convergence properties.

\subsection{Outer convergence}
\autoref{thm:sqpConv} establishes convergence if~\eqref{eq:qpk} is solved exactly in Step~\ref{stp:7} of Algorithm~\ref{alg:basicSQP}. 
Instead, we allow inexact solutions of~\eqref{eq:qpk} and rely on the inexact Newton stopping criterion \cite{Dembo1982}\vspace*{-0.1cm}
\begin{equation} \label{eq:inexactNewton}
\| F^k + \nabla F^{k\top} d^k \| \leq \eta^k \| F^k \|\vspace*{-0.2cm}
\end{equation} with $0 < \eta^k < 1$.
\begin{lem}[Local convergence with inexact SQP steps]\label{lem:inexactNewton}
	Let Assumption~\ref{ass:kkt} hold. 
	Solve~\eqref{eq:qpk} inexactly in Step~\ref{stp:7} of Algorithm~\ref{alg:basicSQP}, form $d^k$ according to~\eqref{eq:d}, and let $d^k$ satisfy~\eqref{eq:inexactNewton} for all $k\geq 0$. 
	Then there exist constants $\varepsilon_2 >0$ and $\eta \in (0,1)$ such that, for all $p^0 \in \mathcal{B}_{\varepsilon_2} ( p^\star)$, the following holds:
	\begin{enumerate}
		\item[i)\hspace*{-0.15cm}] If $\eta^k \leq \eta$ for all $k \geq 0$, then the sequence $\{p^k\}$ generated by Algorithm~\ref{alg:basicSQP} converges q-linearly to~$p^\star$. 
		\item[ii)\hspace*{-0.15cm}] If additionally $\eta^k \rightarrow 0$, then the convergence rate is q-superlinear.
		\item[iii)\hspace*{-0.15cm}] If additionally $\eta^k =O(\| F^k \|)$, then  the convergence rate  is  q-quadratic.\hfill $\square$
	\end{enumerate}
\end{lem} 
\begin{proof}
	We first choose $\varepsilon_2 \in (0,\varepsilon_1]$ with $\varepsilon_1$ from Theorem~\ref{thm:sqpConv}. Then, 
	$\nabla F^\top$ is regular inside $\mathcal{B}_{\varepsilon_2} (p^\star)$.
	To obtain statement~i), we modify~\cite[Theorem 2.3]{Dembo1982} which proves the following: For any $t \in (\eta,1)$, there exists $\varepsilon_2 > 0$ such that, if $p^0 \in \mathcal{B}_{\varepsilon_2} (p^\star)$, then ${\| \nabla F(p^\star)^\top(p^{k+1} - p^\star) \| \leq t \| \nabla F(p^\star)^\top (p^k - p^\star) \|}$ for all $k\geq 0$.
	From the second-to-last inequality follows ${\| p^{k+1} - p^\star \| \leq \text{cond}(\nabla F(p^\star)^\top) t \| p^{k} - p^\star\|}$ \cite[Eq. 7]{Morini99}. 
	Since $\nabla F(p^\star)^\top$ is regular, $\text{cond}(\nabla F(p^\star)^\top)$ is finite and choosing $\eta$, $t$, and $\varepsilon_2$ sufficiently small yields q-linear convergence with ${c = \text{cond}(\nabla F(p^\star)^\top) t < 1}$.		
	Statements ii) and iii) then follow from \cite[Corollary 3.5]{Dembo1982}. 
\end{proof}

The stopping criterion \eqref{eq:inexactNewton} guarantees local convergence of Algorithm~\ref{alg:basicSQP} for iterates starting in $\mathcal B_{\varepsilon_2}$.
However, in practice it is often not known whether $p^0$ lies inside $\mathcal B_{\varepsilon_2}$. 
For an initialization outside  $\mathcal B_{\varepsilon_2}$, it may not be possible to evaluate~\eqref{eq:inexactNewton} as the derivatives of the block rows
$\text{min}(-h_i(x_i^k), \mu_i^k), i \in \mathcal{S}$ in \eqref{eq:F} are ill-defined.
Hence, we develop a modified stopping criterion that is equivalent to~\eqref{eq:inexactNewton} inside $\mathcal{B}_{\varepsilon_2}$, but which can also be evaluated outside $\mathcal{B}_{\varepsilon_2}$.	
We propose the modified stopping criterion\vspace*{-0.1cm}
\begin{equation}\label{eq:modifiedStop}
\| \tilde{F}^k + \nabla \tilde{F}^{k\top} d^k \| \leq \eta^k \| \tilde{F}^k \|,
\end{equation} with

\begin{equation*}
\tilde{F}^k \doteq \begin{bmatrix}
\nabla_{x_1} L_1(x_1^k,\nu_1^k,\mu_1^k,\lambda^k)\\
g_1^k\\
\vdots\\
\nabla_{x_S} L_S(x_S^k,\nu_S^k,\mu_S^k,\lambda^k)\\
g_S^k\\
\left(\sum_{i \in \mathcal{S}}E_i x_i^k \right)- c
\end{bmatrix}
\end{equation*} and where $\nabla \tilde F^{k\top}$ is defined in~\eqref{eq:nablaF}.

\begin{figure*}
\begin{equation} \nabla \tilde F^{k\top} = \begin{bmatrix} 
\nabla_{x_1x_1}^2 L_1^k & \dots  & 0 & \nabla g_1^k &  \dots & 0 & \nabla h_1^k & \dots & 0 & E_1^\top \\
\nabla g_1^{k\top} & \dots & 0 & 0  & \dots & 0 & 0 & \dots & 0 & 0  \\
\vdots & \ddots & \vdots & \vdots & \ddots & \vdots & \vdots & \ddots & \vdots & \vdots\\
0 & \dots & \nabla_{x_Sx_S}^2 L_S^k & 0 & \dots & \nabla g_S^k & 0 & \dots & \nabla h_S^k & E_S^\top\\
0 & \dots & \nabla g_S^{k\top} & 0 & \dots & 0 & 0 & \dots & 0 & 0\\
E_1 & \dots & E_S & 0 & \dots & 0 & 0 & \dots & 0 & 0 
\end{bmatrix} \label{eq:nablaF}
\end{equation}
\hrule 
\end{figure*}
 
Observe that 
$\tilde{F}$ does not include the block rows $\text{min}(-h_i(x_i^k), \mu_i^k), i \in \mathcal{S}$ to avoid differentiability issues outside $\mathcal{B}_{\varepsilon_1}$. 
The subsystems can evaluate~\eqref{eq:modifiedStop} individually and only communicate convergence flags, if $\sum_i E_i (x_i^k + s_i^k) = c$ and if~$\| \cdot \|_\infty$ is chosen.

\begin{lem}[Modified stopping criterion]\label{lem:modStop}
	Suppose Assumption~\ref{ass:kkt} holds, let ${p^k \in \mathcal{B}_{\varepsilon_1}(p^\star)}$, and let $(\bar s^k,\nu^{k+1},\mu^{k+1},\lambda^{k+1})$ be an inexact solution to QP~\eqref{eq:qpk} which has the same active set as $p^\star$, i.e.
	\begin{subequations}\label{eq:correctAS}
		\begin{align}
		[h_i^k]_j + [\nabla h_i^{k\top} \bar s_i^{k}]_j &= 0, \; \forall j \in \mathcal{A}_i, \forall i \in \mathcal{S}\\
		[\mu_i^{k+1}]_j &= 0, \; \forall j \in \mathcal{I}_i,\, \forall i \in \mathcal{S} .
		\end{align}
	\end{subequations} Moreover, set $d^k = (\bar s^k, \Delta \nu^k, \Delta \mu^k, \Delta \lambda^k)$. Then,	
	\[\| \tilde{F}^k  +  \nabla \tilde{F}^{k\top} d^k\| \hspace*{-1mm} \leq \hspace*{-1mm} \eta^k \| \tilde{F}^k \| \hspace*{-2mm} \implies \hspace*{-2mm} \| F^k  +  \nabla F^{k\top} d^k \| \hspace*{-1mm} \leq \hspace*{-1mm} \eta^k \| F^k \|.\]
	\phantom{m}\hfill $\square$	
\end{lem}
\begin{proof}
	From \eqref{eq:correctAS} follows $\| \tilde{F}^k  +  \nabla \tilde{F}^{k\top} d^k\|  = \| F^k + \nabla F^{k\top} d^k\|$. Since $ \| \tilde{F}^k \| \hspace*{-0.1cm} \leq \hspace*{-0.1cm}\| F^k \|$, we obtain the implication.
\end{proof}
		
\subsection{Inner convergence}	
To apply the stopping criterion \eqref{eq:modifiedStop}, we show that ADMM stops at the correct active set in a neighborhood of $p^\star$. 
We first show that the optimal active sets of $\eqref{eq:sepForm}$, $\eqref{eq:reformulatedQP}$, and $\eqref{eq:localQP}$ are equivalent in a neighborhood of $p^\star$, if ADMM is initialized appropriately.
Let $q \doteq (\bar{s},\gamma)$ and $q_i \doteq (\bar{s}_i,\gamma_i)$ for all~$i \in \mathcal{S}$.
\begin{lem}[Active set of the ADMM subsystem QPs]\label{lem:admmAS}
	Suppose Assumption~\ref{ass:kkt} holds and denote the KKT point of the two-block QP~\eqref{eq:reformulatedQP} formulated at $p^k$ by  $(s^{k,\star},\bar{s}^{k,\star},\nu^{k,\star},\mu^{k,\star},\gamma^{k,\star},\lambda^{k,\star})$.
	Then, there exists a constant $\varepsilon_4 > 0$ such that the following holds.
	If $p^k \in \mathcal{B}_{\varepsilon_1}(p^\star)$ and ${q^l \in \mathcal{B}_{\varepsilon_4}(q^{k,\star})}$, then the solution $(s_i^{l+1},\mu_i^{l+1})$ to the subsystem QP~\eqref{eq:localQP} with parameters $q_i^l$ has the same active set as $p^\star$, i.e.
	\begin{align*}
	\hspace{0.85cm}[h_i^k]_j + [\nabla h_i^{k\top} {s}_i^{l+1}]_j &= 0, \;\; \forall j \in \mathcal{A}_i,\; \forall i \in \mathcal{S},\\
	[{\mu}_i^{l+1}]_j &= 0, \;\; \forall j \in \mathcal{I}_i,\; \; \forall i \in \mathcal{S}. \hspace*{0.9cm} \square
	\end{align*} 
\end{lem}
\begin{proof}
	We first show that, (a), for $p^k \approx p^\star$, the solution of the ADMM subsystem QP~\eqref{eq:localQP} has the correct active set, if \eqref{eq:localQP} is parameterized with $q^{k,\star}$. We then show that, (b), this also holds for small variations in $q$.
		
	(a) As derived in part (a) of the proof of Theorem~\ref{thm:sqpConv}, the solution $p^{k,\star}$ to QP~\eqref{eq:qpk} is unique and has the same active set as $p^\star$ for all $p^k \in \mathcal{B}_{\varepsilon_1}(p^\star)$.
	Comparing the KKT conditions of QPs~\eqref{eq:qpk} and~\eqref{eq:reformulatedQP}, we see that if $(s^{k,\star},\nu^{k,\star},\mu^{k,\star},\lambda^{k,\star})$ is the corresponding KKT point of QP~\eqref{eq:qpk}, then $(s^{k,\star},\bar{s}^{k,\star},\nu^{k,\star},\mu^{k,\star},\gamma^{k,\star},\lambda^{k,\star})$ with ${\bar{s}^{k,\star} = s^{k,\star}}$ and ${\gamma^{k,\star} = E^\top \lambda^{k,\star}}$ with ${E \doteq \begin{bmatrix} E_1 & \hspace*{-0.1cm}\dots\hspace*{-0.1cm} & E_S \end{bmatrix}}$ is the KKT point of QP~\eqref{eq:reformulatedQP}.
	Moreover, if $(s^{k,\star},\bar{s}^{k,\star},\nu^{k,\star},\mu^{k,\star},\gamma^{k,\star},\lambda^{k,\star})$ is the KKT point of \eqref{eq:reformulatedQP}, then $(s_i^{k,\star},\nu_i^{k,\star},\mu_i^{k,\star})$ is the KKT point of QP~\eqref{eq:localQP} for subsystem $i \in \mathcal{S}$ with parameters $q_i^{k,\star}$.
	Thus, for all $i \in \mathcal{S}$, if QP~\eqref{eq:localQP} is parameterized by $q_i^{k,\star}$, then $(s_i^{k,\star},\nu_i^{k,\star},\mu_i^{k,\star})$ is a unique solution with active set $\mathcal{A}_i$. 
 
	(b) From (a), we have that the KKT point $(s_i^{k,\star},\nu_i^{k,\star},\mu_i^{k,\star})$ of QP~\eqref{eq:localQP} for subsystem $i \in \mathcal{S}$ with parameters $q_i^{k,\star}$ satisfies LICQ and strict complementarity.
	Moreover, the Hessians $H_i^k$ satisfy~\eqref{eq:sconvQP} for all $i \in \mathcal{S}$.
	Hence, we can invoke the BST and obtain that solving the subsystem QP returns the correct active set if $q^l \in \mathcal{B}_{\varepsilon_4}(q^{k,\star})$ for some $\varepsilon_4 > 0$.
	That is, if $p^k \in \mathcal{B}_{\varepsilon_1}(p^\star)$ and if $q^l \in \mathcal{B}_{\varepsilon_4}(q^{k,\star})$, then	
	\begin{align*}
	\hspace{0.85cm}[h_i^k]_j + [\nabla h_i^{k\top} {s}_i^{l+1}]_j &= 0, \;\; \forall j \in \mathcal{A}_i,\; \forall i \in \mathcal{S},\\
	[{\mu}_i^{l+1}]_j &= 0, \;\; \forall j \in \mathcal{I}_i,\; \; \forall i \in \mathcal{S}.
	\end{align*}
\end{proof}
Lemma~\ref{lem:admmAS} shows that the active set of the iterate $s^{l+1}$ found by solving the subsystem QPs is correct if $p^k \approx p^\star$ and $q^l \approx q^{k,\star}$.
We next show that this property carries over to the iterate $\bar{s}^{l+1}$ produced by the averaging Step~\ref{eq:ADMMStep2} of ADMM.
To this end, we make the following assumption to ensure that the averaging step does not affect variables which are bound by active inequalities.
Let the map ${h(x) \doteq (h_1(x_1),\dots,h_S(x_S))}$ denote the centralized inequality constraints and define the ADMM averaging matrix 
\[M \doteq (I - E^\top (E E^\top)^{-1} E).\]

\begin{ass}[Decoupled inequality constraints]\label{ass:decoupled}
	The inequality and coupling constraints satisfy $\nabla h(x)^\top E^\top = 0$ for all $x \in \mathbb{R}^n$. Furthermore, the right-hand side of the coupling constraint~\eqref{eq:consConstr} is $c = 0$. \hfill $\square$
\end{ass}

Assumption~\ref{ass:decoupled} states that decision variables which participate in inequality constraints are not directly coupled to neighbors.
While this assumption facilitates the below convergence proof, d-SQP can also be applied without convergence guarantees if Assumption~\ref{ass:decoupled} does not hold.
We further note that inequality constraints can always be decoupled by introducing additional decision variables as shown in Example~\ref{ex:decoupling} in the appendix.

\begin{lem}[Active set of the ADMM averaging step]\label{lem:admmConvergence}
	Let Assumptions~\ref{ass:kkt} and \ref{ass:decoupled} hold and let $p^k \in \mathcal{B}_{\varepsilon_1}(p^\star)$.
	Furthermore, let the ADMM dual initialization satisfy $M \gamma^0 = 0$.
	Then, there exists a constant $\varepsilon_3 > 0$ such that if ${q^0 \in \mathcal{B}_{\varepsilon_3}(q^{k,\star})}$, 
	then the iterates $(\bar{s}^{l+1},\mu^{l+1})$ produced by Algorithm~\ref{alg:ADMM} have the same active set as $p^\star$ for all $l\geq 0$, i.e.
	\begin{align*}
	\hspace*{0.7cm}[h_i^k]_j + [\nabla h_i^{k\top} \bar{s}_i^{l+1}]_j &= 0, \forall j \in \mathcal{A}_i, \forall i \in \mathcal{S}, \forall l\geq 0,\\
	\qquad \quad [{\mu}_i^{l+1}]_j &= 0, \forall j \in \mathcal{I}_i, \; \forall i \in \mathcal{S}, \forall l\geq 0. \hspace*{0.25cm} \square	
	\end{align*} 
\end{lem}
\begin{proof}
	We first show that, (a), inside $\mathcal{B}_{\varepsilon_1}(p^\star)$ problem~\eqref{eq:reformulatedQP} is convex on the constraints. 
	Then, (b), standard ADMM convergence results imply that the iterates will stay in a bounded neighborhood of $q^{k,\star}$.
	Finally, (c), we show that this neighborhood can be chosen such that the ADMM iterates have the correct active set.

	(a) Let $G_i^k \doteq \nabla g_i^{k\top}$. Recall that for $A \in \mathbb{R}^{m \times n}$ the null space of $A$ and the range space of $A^\top$ together form $\mathbb{R}^n$ \cite[p. 603]{Nocedal2006}. We can therefore write ${s_i = u_i + v_i}$, where $u_i$ lies in the null space of $G_i^k$ and $v_i$ lies in the range space of $G_i^{k\top}$.
	From \eqref{eq:locEqCon} and LICQ follows ${v_i = - G_i^{k\top} (G_i^{k} G_i^{k\top})^{-1} g_i^k}$.  
	Since $v_i$ is uniquely determined and since ${u_i^\top H_i^k u_i > 0}$ because of Assumption~\ref{ass:kkt}, the local objective 
	$\phi_i^k(s_i)$ is a convex function for all feasible~$s_i$. 	
	
	(b) ADMM is therefore guaranteed to converge to ${(s^{k,\star},\bar{s}^{k,\star},\nu^{k,\star},\mu^{k,\star},\gamma^{k,\star},\lambda^{k,\star})}$ and the ADMM iterates satisfy \cite[Eq. (3.5)]{He2015}
	\begin{align*}\label{eq:ADMMconvergence}
	\rho \| \bar{s}^{l+1}-\bar{s}^{k,\star} \|_2^2 + &\frac{1}{\rho} \| \gamma^{l+1}-\gamma^{k,\star} \|_2^2 \leq\\
	&\rho \| \bar{s}^{l} - \bar{s}^{k,\star} \|_2^2 + \frac{1}{\rho} \| \gamma^{l}-\gamma^{k,\star} \|_2^2
	\end{align*} for all $l\geq0$. 
	Since $0 < \rho < \infty$, the above inequality implies that, for all $\varepsilon > 0$, there exists $\varepsilon_3 > 0$ such that if $q^0 \in \mathcal{B}_{\varepsilon_3}(q^{k,\star})$, then $q^{l} \in \mathcal{B}_{\varepsilon}(q^{k,\star})$ for all $l \geq 0$.
	
	(c) From (b), we can choose $\varepsilon_3 > 0$ such that $q^{l} \in \mathcal{B}_{\varepsilon_4}(q^{k,\star})$ for all $l \geq 0$ with $\varepsilon_4$ from Lemma~\ref{lem:admmAS}.
	Thus, the iterates $(s_i^{l+1},\mu_i^{l+1})$ remain at the correct active set for all ADMM steps, i.e.
	\begin{align*}
	\hspace*{0.7cm}[h_i^k]_j + [\nabla h_i^{k\top} {s}_i^{l+1}]_j &= 0, \forall j \in \mathcal{A}_i, \forall i \in \mathcal{S}, \forall l\geq 0,\\
	\qquad \quad [{\mu}_i^{l+1}]_j &= 0, \forall j \in \mathcal{I}_i, \; \forall i \in \mathcal{S}, \forall l\geq 0. 
	\end{align*}
	The KKT system of the equality constrained QP in the averaging Step~\ref{eq:ADMMStep2} yields~\cite{Stomberg2022}
	\[
	\bar{s}^{l+1} = Ms^{l+1} = (I - E^\top (E E^\top)^{-1} E) s^{l+1}.
	\]
	Hence, we obtain
	\begin{align*}
	h^k + &\nabla h^{k\top}\bar{s}^{l+1} = h^k + \nabla h^{k\top}(s^{l+1} + \bar{s}^{l+1} - s^{l+1})\\
	& = h^k + \nabla h^{k\top}(s^{l+1} + M s^{l+1} - s^{l+1})\\
	& = h^k + \nabla h^{k\top}s^{l+1} - \underbrace{\nabla h^{k\top} E^\top}_{0,\text{ Ass.~\ref{ass:decoupled}}} (E E^\top)^{-1} E s^{l+1}\\
	&= h^k + \nabla h^{k\top} s^{l+1}.
	\end{align*} Inserting the correct active set for $s^{l+1}$ yields
	\begin{align*}
	\hspace*{0.7cm}[h_i^k]_j + [\nabla h_i^{k\top} \bar{s}_i^{l+1}]_j &= 0, \forall j \in \mathcal{A}_i, \forall i \in \mathcal{S}, \forall l\geq 0.
	\end{align*} 
	The averaged variable $\bar{s}^{l+1}$ hence lies at the correct active set.	
\end{proof}

\begin{rem}[ADMM dual initialization]\label{rem:dualinit}
Lemma~\ref{lem:admmConvergence} states that the ADMM initialization satisfies $M\gamma^0 = 0$. 
This is the case if $\gamma_i^\top = E_i^\top \lambda^0$ for all $i \in \mathcal{S}$, because
\[
M \gamma^0 = M E^\top \lambda^0 = (E^\top - E^\top) \lambda^0 = 0 \quad \text{for all} \quad \lambda^0 \in \mathbb{R}^{n_c}.
\]
Moreover, the ADMM updates ensure that~\cite[Ch. 7]{Boyd2011}
\[
M \gamma^l = 0 \implies M \gamma^{l+1} = 0.
\] 
\end{rem} \vspace*{-0.4cm} \hfill $\square$

\subsection{Local convergence of decentralized SQP}
\setlength{\textfloatsep}{14pt}
\begin{algorithm}[t]
	\caption{d-SQP for solving \eqref{eq:sepForm}}
	\begin{algorithmic}[1]
		\State SQP initialization: $k=0,\lambda^0$, $(x_i^0,\nu_i^0,\mu_i^0,\gamma_i^0 = E_i^\top \lambda^0)$ for all $i \in \mathcal{S}$, $\eta^0 < 1,\epsilon$ \label{dsqp-stp:1}
		\While{ $\| F^k \| \nleq \epsilon$} \label{dsqp-step:2}
		\State compute $\nabla f_i^k, g_i^k, \nabla g_i^k, h_i^k, \nabla h_i^k, H_i^k$ for all $i \in \mathcal{S}$ \label{stp:3}
		\State compute $\tilde{F}_i^k$ and $\nabla \tilde{F}_i^k$ for all $i \in \mathcal{S}$\label{dsqp-step4}
		\State ADMM initialization: $l=0$, $(\bar{s}_i^l, \gamma_i^l) = (0,\gamma_i^k)$, and
		
		\hspace*{-3mm} $(\nu_i^l, \mu_i^l) = (\nu_i^k,\mu_i^k)$ for all $i \in \mathcal{S}$
		 \While{ $l = 0$ \textbf{or} $\| \tilde{F}^k + \nabla \tilde{F}^{k\top} d^l \| \nleq \eta^k \| \tilde{F}^k \|$}\label{dsqp-step:6}		
		\State \hspace*{-3mm} $\displaystyle (s_i^{l+1},\nu_i^{l+1},\mu_i^{l+1})  \leftarrow  \min_{s_i \in \mathbb{S}_i^k} L_{\rho,i}^k(s_i,\bar{s}_i^l,\gamma_i^l) \, \forall i \in \mathcal{S}$  
		\State \hspace*{-3mm} $\displaystyle \bar{s}^{l+1} = \argmin_{\bar{s} \in \mathbb{E}} \sum_{i \in \mathcal{S}} L_{\rho,i}^k(s_i^{l+1},\bar{s}_i,\gamma_i^{l})$ \label{eq:sqp_comm}
		\State \hspace*{-3mm} $\gamma_i^{l+1} = \gamma_i^l + \rho (s_i^{l+1}-\bar{s}_i^{l+1})$ for all $i \in \mathcal{S}$
		\State \hspace*{-3mm} $l \leftarrow l+1$
		\EndWhile
		\State $(x_i^{k+1}\hspace*{-1mm},\nu_i^{k+1}\hspace*{-1mm},\mu_i^{k+1}\hspace*{-1mm},\gamma_i^{k+1}) \hspace*{-1mm} = \hspace*{-1mm} (x_i^k + \bar{s}_i^l,\nu_i^{l},\mu_i^{l},\gamma_i^l)$  $\forall i \in \mathcal{S}$
		\State choose $\eta^{k+1} \leq \eta^k$
		\State $k \leftarrow k+1$
		\EndWhile\label{euclidendwhile}
		\State \textbf{return} $x_i^k$ for all $i \in \mathcal{S}$

	\end{algorithmic} \label{alg:d-SQP}

\end{algorithm}

Algorithm~\ref{alg:d-SQP} summarizes the decentralized SQP method~(d-SQP). The algorithm follows a bi-level structure and we denote the SQP (outer) iterations by~$k$ and the ADMM (inner) iterations by $l$.
Within each SQP iteration, ADMM is terminated via the inexact Newton-type stopping criterion in Step~\ref{dsqp-step:6} based on the SQP step $d^l \doteq (\bar s^l, \nu^l - \nu^k, \mu^l - \mu^k, \lambda^l - \lambda^k)$ and we will comment on the evaluation of $\lambda^l - \lambda^k$ in implementations below. 
We are now ready to state our main result.
\begin{thm}[Local convergence of d-SQP]\label{thm:dsqpConv}
	Let Assumptions~\ref{ass:kkt} and \ref{ass:decoupled} hold. 
	Then, there exist constants $\varepsilon > 0$ and $\eta > 0$ such that, for all $p^0 \in \mathcal{B}_{\varepsilon} ( p^\star)$, the following holds:
	\begin{enumerate}
		\item[i)\hspace*{-0.15cm}] If $\eta^k \leq \eta$ for all $k \geq 0$, then the sequence $\{p^k\}$ generated by Algorithm~\ref{alg:d-SQP} converges to~$p^\star$ and the convergence rate is q-linear in the outer iterations.
		\item[ii)\hspace*{-0.15cm}] If additionally $\eta^k \rightarrow 0$, then the convergence rate is q-superlinear in the outer iterations.
		\item[iii)\hspace*{-0.15cm}] If additionally $\eta^k = O(\| \tilde{F}^k \|)$, then the convergence rate is q-quadratic in the outer iterations. \hfill$\square$
		\addtolength{\itemindent}{-0.4cm}
	\end{enumerate}
\end{thm}
\begin{proof}
	We first show that, (a), the ADMM initialization $(\bar s^0,\gamma^0) = (0,E^\top\lambda^k)$ lies within a neighborhood of the subproblem solution $q^{k,\star}$ for all SQP iterations. We then show that, (b), this neighborhood can be chosen such that the ADMM iterations are at the correct active set at all iterations. Finally, (c), we invoke Lemmas~\ref{lem:inexactNewton}--\ref{lem:modStop} to prove convergence.
	
	(a) We first choose $\varepsilon$ such that $0 < \varepsilon \leq \varepsilon_2$. The convergence of the Newton and inexact Newton methods implies that if ${p^k \in \mathcal{B}_{\varepsilon}(p^\star)}$, then $p^{k+1} \in \mathcal{B}_{\varepsilon}(p^\star)$ holds for the Newton method as well as for the inexact Newton method~{\cite[Theorem 2.3]{Dembo1982}}.
	Therefore $\| p^{k+1} - p^k \| \leq 2 \varepsilon.$
	Recall that ${\| (x,y) \| \geq \| x \| }$ for any vectors $x$ and $y$.
	Hence, ${\|(x^{k+1}-x^k,\lambda^{k+1}-\lambda^k)\| \leq 2\varepsilon}$.
	Further recall that ${\gamma^k = E^\top \lambda^k}$.
	Hence, $\| (x^{k+1}-x^k,\gamma^{k+1}-\gamma^k ) \| \leq \varepsilon_{5}$ with $\varepsilon_{5} = \text{max}(\| E^\top \|,1)\cdot 2 \varepsilon$. 	
	That is, ${(0,\gamma^k) \in \mathcal{B}_{\varepsilon_{5}}((x^{k+1}-x^k, \gamma^{k+1}))}$.
	The iterations of the Newton method read ${x^{k+1}= x^k+ \bar{s}^{k,\star}}$ and ${\lambda^{k+1}\hspace*{-0.1cm}=\hspace*{-0.1cm}\lambda^{k,\star}}$ and we hence get ${(0,E^\top\lambda^k) \in \mathcal{B}_{\varepsilon_{5}}(q^{k,\star})}$.
	 
	(b) The d-SQP initialization $\gamma_i^\top = E_i^\top \lambda^0$ for all $i \in \mathcal{S}$ yields $M \gamma^0 = 0$ for all $\lambda \in \mathbb{R}^{n_c}$.
	Moreover, the ADMM updates ensure $M \gamma^l = 0$ for all $l \geq 0$~\cite[Ch. 7]{Boyd2011}. 
	Lemma~\ref{lem:admmConvergence} hence shows that the ADMM iterates $(\bar{s}^{l+1},\mu^{l+1})$ are at the correct active set for all $l\geq0$, if $q^0 \in \mathcal{B}_{\varepsilon_3}(q^{k,\star})$. 	
	We therefore choose ${\varepsilon  \leq \text{min}(\varepsilon_2,\varepsilon_3 / (2 \text{max}(1,\|E^\top \|)))}$ such that $\varepsilon_5 \leq \varepsilon_3$.
	The ADMM initialization then satisfies ${q^0 = (0,E^\top \lambda) \in \mathcal{B}_{\varepsilon_3}(q^{k,\star})}$.
	By Lemma~\ref{lem:admmConvergence}, the iterates $(\bar{s}^{l+1},\mu^{l+1})$ therefore are at the correct active set for all ADMM iterations $l \geq 0$ and for all SQP iterations $k \geq 0$.
	 
	(c) The ADMM convergence invoked in the proof of Lemma~\ref{lem:admmConvergence} ensures that ADMM can satisfy the stopping criterion \eqref{eq:modifiedStop}. 
	By Lemma~\ref{lem:modStop}, satisfaction of~\eqref{eq:modifiedStop} together with the correct active set implies satisfaction of the inexact Newton stopping criterion~\eqref{eq:inexactNewton}.
	Lemma~\ref{lem:inexactNewton} then implies local q-convergence in the outer iterations.
	Recall that $\| \tilde F^k \| \leq \| F^k\|$. Hence $\eta^k = O(\|\tilde F^k \|) \implies \eta^k = O(\| F^k \|)$ such that we obtain the q-quadratic convergence in statement iii$)$. This concludes the proof.
\end{proof}

\subsection{Communication requirements and discussion}

Four steps of Algorithm~\ref{alg:d-SQP} require communication between subsystems: Step~\ref{dsqp-stp:1} to initialize $\lambda^0$, Steps~\ref{dsqp-step:2} and \ref{dsqp-step:6} to evaluate the stopping criteria, and Step~\ref{eq:sqp_comm} for the $\bar{s}$ update in ADMM.

It is well known that Step~\ref{eq:sqp_comm} may be computed efficiently as a decentralized averaging step, if NLP~\eqref{eq:sepForm} is a consensus-type problem, cf.~\cite{Boyd2011,Stomberg2022}. 
That is, the update of $\bar{s}$ in Step~\ref{eq:sqp_comm} then requires only the communication of vectors between neighboring subsystems if we choose  ${\gamma_i^0 = E_i^\top \lambda^0}$. 

If $\| \cdot \|_\infty$ is chosen, then the stopping criteria can be evaluated locally and in Steps~\ref{dsqp-step:2} and \ref{dsqp-step:6} only convergence flags must be communicated, because the update $x^{k+1} = x^k + \bar{s}^k$ ensures $\sum_{i \in \mathcal{S}} E_i x_i^k = c$ for all $k \geq 1$.

\begin{rem}[Dual iterates of the coupling constraints]
	Theorem~\ref{thm:dsqpConv} proves convergence with the dual variable $\lambda^k$.
	However, instead of computing $\lambda^k$ explicitly, Algorithm~\ref{alg:d-SQP} exploits that $E_i^\top \lambda = \gamma_i$ to evaluate Steps~\ref{dsqp-step:2}, \ref{dsqp-step4}, and~\ref{dsqp-step:6}. \hfill $\square$
\end{rem}

\begin{rem}[Hessian regularization]
	ADMM is guaranteed to converge if the cost functions $\phi_i^k(s_i)$ are convex on the constraints~$\mathbb{S}_i^k$. 
	That is, we do not require $H_i^k$ to be positive definite, but only positive definite on the equality constraint null space. 
	However, this may not be the case if d-SQP is initialized far away from $p^\star$. To this end, we regularize the reduced Hessian via \cite[Eq. 3.43]{Nocedal2006} with ${\delta = 10^{-4}}$, where $\delta$ is in the notation of \cite{Nocedal2006}. \hfill $\square$
\end{rem}

\begin{rem}[Inner method and application to QPs]
Here, we present d-SQP with ADMM as an inner method for solving~\eqref{eq:qpk}. 
Instead of ADMM, other decentralized methods with guaranteed convergence to a KKT point for convex QPs can also be used. If~\eqref{eq:sepForm} itself is a QP, then d-SQP is equivalent to applying the inner method to a possibly regularized version of~\eqref{eq:sepForm}. $\hfill \square$
\end{rem}

\section{Numerical Results}\label{sec:numerics}

We compare the performance of d-SQP to four other methods for an Alternating Current OPF problem for the IEEE 118-bus system.
The first method for comparison is standalone ADMM, which solves~\eqref{eq:sepForm} directly and which we denote as ADMM in the following.
The second and third methods are bi-level ALADIN variants, where the ALADIN coordination~QP is solved via essentially decentralized Conjugate Gradients (d-CG) or via an ADMM variant (d-ADMM).
The fourth method is an essentially decentralized Interior Point method~(d-IP)~\cite{Engelmann2021a}.

The problem consists of four subsystems with  a total of $n \hspace*{-0.1cm}=\hspace*{-0.1cm} 576$ variables, $n_g \hspace*{-0.1cm}=\hspace*{-0.1cm} 470$ non-linear equality constraints, ${n_h\hspace*{-0.1cm} = \hspace*{-0.1cm}792}$ linear inequality constraints, and $n_c\hspace*{-0.1cm}=\hspace*{-0.1cm}52$ coupling constraints. 
We initialize voltage magnitudes as $1$ and remaining variables as $0$. 

For ADMM, we tune $\rho$ for fastest convergence in the set $\{100,700,800,900,10^3,10^4\}$ and obtain $\rho = 800$. 
The respective subproblems in ADMM and ALADIN are solved with \texttt{IPOPT}~\cite{Biegler2009} and the QPs in d-SQP are solved with \texttt{qpOASES}~\cite{Ferreau2014}.
For d-SQP, we choose $\eta^0 = 0.8$, ${\eta^{k+1} = 0.9 \eta^k}$, select $\rho = 700$ from the set $\{600,700,800\}$, and initialize the inner iterations with $(\bar{s}_i^0,\gamma_i^0) = (0,E_i^\top \lambda^k)$.
For bi-level ALADIN, we run either 70 d-CG iterations or 100 d-ADMM iterations per ALADIN iteration, tune $\rho$ in the set $\{10,50,100,250,200,1000\}$, and obtain $\rho = 150$.
For d-IP, we use the parameters from \cite{Engelmann2021a}.
We use \texttt{CasADi}~\cite{Andersson2019} to compute derivatives via algorithmic differentiation in all methods.

\begin{figure}
\includegraphics[width=\columnwidth]{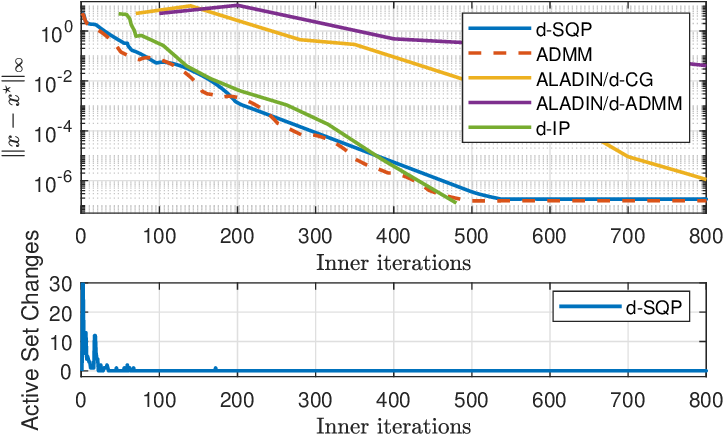}
\vspace*{-0.5cm}\caption{d-SQP, ADMM, bi-level ALADIN, and d-IP for 118-bus OPF.}\label{fig:x118bus}
\end{figure}

The top part of Figure~\ref{fig:x118bus} shows the convergence to the minimizer found by \texttt{IPOPT}.\footnote{The convergence results for ADMM shown in Figure~\ref{fig:x118bus} differ from the results shown in \cite{Engelmann2021a}. Here we use a different formulation of the augmented Lagrangian to facilitate the averaging step, cf. \cite{Stomberg2022}.} For 
d-SQP, bi-level ALADIN, and ${\text{d-IP}}$ we count the inner iterations.
Figure~\ref{fig:x118bus} shows that d-SQP requires more iterations than ADMM to achieve an equivalent accuracy in $\|x-x^\star\|_\infty$.
However, at the level of each subsystem, ${\text{d-SQP}}$ only solves QPs whereas ADMM solves NLPs.
As a result, in our prototypical implementation ADMM takes $48\,$s to achieve ${\|x-x^\star\|_\infty<10^{-6}}$, whereas d-SQP only takes $11\,$s.
Moreover, d-SQP is guaranteed to converge locally for this problem, whereas ADMM might diverge~\cite{Christakou2017}. 
The bottom plot of Figure~\ref{fig:x118bus} shows the number of inequality constraints which toggle between being active and inactive.
The active set settles within 200 iterations, which indicates that the area of local convergence is reached. 

All methods require local neighbor-to-neighbor communication of the same complexity: $2n_c$ floats for all subsystems combined per inner iteration. 
ALADIN/d-CG and d-IP further communicate two floats globally per d-CG iteration. 
The faster convergence of d-IP in Figure~\ref{fig:x118bus} therefore comes at the cost of global communication with low complexity.
Hence, the bi-level SQP scheme shows competitive \textit{decentralized} performance \textit{and} it admits local convergence guarantees. 

\section{Conclusion and Outlook}

This paper has established new convergence guarantees for decentralized SQP schemes under non-convex constraints.
The proposed method solves quadratic subproblems with ADMM.
In contrast to existing decentralized SQP methods, we allow for inexact ADMM solutions via an appropriate decentralized stopping criterion for the inner iterations. 
Numerical results show competitive performance to further optimization schemes for an example from power systems.
Future work will consider the globalization of the method and its application to distributed model predictive control.

	\renewcommand*{\bibfont}{\footnotesize}
	\footnotesize
	\printbibliography
	\normalsize	
	
\appendix

\begin{ex}[Reformulation of coupled inequalities]\label{ex:decoupling}
	Consider the partially separable QP with two subsystems,
	\begin{subequations}\label{eq:ex}
	\begin{align}
	\min_{x_1 \in \mathbb{R},x_2 \in \mathbb{R}} 10(x_1 - 10)^2 &+ (x_2 - 1)^2\\
	\text{subject to} \quad x_1 - 1 &\leq 0 \; | \; \mu_1,\\
	x_1 - x_2 &= 0 \; | \; \lambda.\label{eq:ex:coup}
	\end{align}
	\end{subequations} Here, $h(x) = x_1 - 1$ and $E = \begin{bmatrix} 1 & -1 \end{bmatrix}$ such that 
	\[
	\nabla h(x)^\top E^\top = \begin{bmatrix} 1 & 0 \end{bmatrix} \cdot \begin{bmatrix}  1 & -1 \end{bmatrix}^\top = 1.
	\]
	That is, Assumption~\ref{ass:decoupled} does not hold, because the variable $x_1$ participates in the inequality constraint \textit{and} the coupling constraint~\eqref{eq:ex:coup}.
	For QP~\eqref{eq:ex}, d-SQP reduces to ADMM and can directly be applied without any problem reformulations. 
	However, the convergence guarantee in Theorem~\ref{thm:dsqpConv} does not apply.
	If one wishes to invoke the convergence theorem, then Assumption~\ref{ass:decoupled} can be fulfilled by adding a decision variable to subsystem one and by rewriting the QP as
	\begin{subequations}\label{eq:ex2}
		\begin{align*}
		\min_{x_1 \in \mathbb{R}^2,x_2 \in \mathbb{R}} 10\left( \begin{bmatrix} 1 & 0 \end{bmatrix} x_1 - 10 \right)^2 &+ (x_2 - 1)^2\\
		\text{subject to} \quad \begin{bmatrix} 1 & -1 \end{bmatrix} x_1 &= 0 \; | \; \nu_1,\\
		\begin{bmatrix} 1 & 0 \end{bmatrix} x_1 - 1 &\leq 0 \; | \; \mu_1,\\
		\begin{bmatrix} 0 & 1 \end{bmatrix} x_1 - x_2 &= 0 \; | \; \lambda.
		\end{align*}
	\end{subequations}
	Then
	\[
	\nabla h(x)^\top E^\top = \begin{bmatrix} 1 & 0 & 0 \end{bmatrix} \cdot \begin{bmatrix} 0 & 1 & -1 \end{bmatrix}^\top = 0
	\] and Assumption~\ref{ass:decoupled} holds. \hfill $\square$	
\end{ex}

\end{document}